\documentclass{ifacconf}

\def\rep#1{(\ref{#1})}

\newcommand{\R}{\mathbb{R}}

\def\send#1#2{\stackrel{#1}{\hbox to #2{\rightarrowfill}}}
\def\-{\!\!\!\!\!-}

 \def\qed{ \rule{.1in}{.1in}}
\def\eq#1{\begin{equation}#1\end{equation}}

\newcommand{\matt}[1]{\begin{bmatrix}#1\end{bmatrix}}

\newcommand{\dfb}{\stackrel{\Delta}{=}}

\def\qed{ \rule{.1in}{.1in}}

\def\R{{\rm I\!R}} 

\newcounter{seqn}[equation]
\def\theseqn{\arabic{equation}\alph{seqn}}

\def\endseqn{\eqno \@seqnnum
$$\ignorespaces}
\def\@seqnnum{(\theseqn)}

\newskip\mcentering \mcentering=0pt plus 1000pt minus 1000pt

\def\meqalignno#1{
\halign to\displaywidth{
    \hbox to 0pt{\kern\displaywidth\llap{$##$}\hss}\tabskip=\mcentering
    &\hfil$\displaystyle{##}$\tabskip=\mcentering
   &&$\displaystyle{{}##}$\hfil\tabskip=\mcentering
    \crcr
    #1\crcr}}

\def\rep#1{(\ref{#1})}

\def\eq#1{\begin{equation}#1\end{equation}}

\def\dspace{\multiply\normalbaselineskip 150
		  \divide\normalbaselineskip 100 \normalbaselines
		  \csname @@normalbaselineskip\endcsname\normalbaselineskip}
\def\sspace{\multiply\normalbaselineskip 200
		 \divide\normalbaselineskip 300 \normalbaselines
		 \csname @@normalbaselineskip\endcsname\normalbaselineskip}
\def\sdspace{\multiply\normalbaselineskip 160
		 \divide\normalbaselineskip 150 \normalbaselines
		 \csname @@normalbaselineskip\endcsname\normalbaselineskip}

\def\@{\tilde}
\def\3dot#1{\buildrel\textstyle...\over#1}

\usepackage{graphicx}      
\usepackage{natbib}        
\usepackage{amsfonts}
\usepackage{amsmath}
\begin{document}
\begin{frontmatter}

\title{A Distributed Algorithm for Computing a Common Fixed Point of a Family of Paracontractions}

\thanks{The problem to which this paper is addressed was prompted by useful discussions with Bronislaw Jakubczyk.
 This work was supported by the US
Air Force Office of Scientific Research and by the National Science Foundation.}

\author[First]{Daniel Fullmer}
\author[First]{Lili Wang}
\author[First]{A. Stephen Morse}

\address[First]{Department of Electrical Engineering,
New Haven, CT 06520 USA (e-mail: \{daniel.fullmer, lili.wang, as.morse\}@yale.edu)}

\begin{abstract}

A distributed algorithm is described for finding a  common fixed point of  a family of $m>1$  nonlinear maps
 $M_i:\R^n\rightarrow\R^n$
assuming that each  map is a paracontraction and that such a common fixed point exists.
 The  common fixed point  is
simultaneously computed  by $m$ agents
assuming each agent $i$ knows only  $M_i$,
 the current estimates  of the fixed point generated by its neighbors, and nothing more. Each agent
recursively updates its estimate  of the fixed point by utilizing the current estimates
  generated by each of  its neighbors.
  Neighbor relations are characterized by a  time-dependent directed graph $\mathbb{N}(t)$ whose
 vertices correspond to agents and whose arcs depict neighbor relations.
It is  shown that for
any family of paracontractions $M_i,\;i\in\{1,2,\ldots ,m\}$ which has at least one common fixed point,
and any sequence of strongly connected neighbor graphs $\mathbb{N}(t)$, $t=1,2,\ldots$, the algorithm causes
all agent estimates to converge to a common fixed point.
\end{abstract}


\end{frontmatter}

\section{Introduction}

This paper is concerned with distributed algorithms for  enabling a group of $m>1$ mobile autonomous
agents to solve certain
types of nonlinear equations over a network. It is assumed that  each  agent can receive information from
its neighbors where by a {\em neighbor} of agent $i$ is meant any
other agent within agent $i$'s reception range.
We write $\mathcal{N}_i(t)$
for the labels of agent $i$'s neighbors at time $t$, and we always
take agent $i$ to be a neighbor of itself. Neighbor relations at
time $t$ can be conveniently characterized by a directed graph
$\mathbb{N}(t)$ with $m$ vertices and a set of arcs defined so that there is
an arc in $\mathbb{N}(t)$ from vertex $j$ to vertex $i$ just in case agent $j$
is a neighbor of agent $i$ at time $t$.  Each agent $i$ has a
real-time dependent state vector $x_i(t)$ taking values in $\R^n$, and
we assume that the information agent $i$ receives from neighbor
$j$ at time $t$ is $x_j(t)$. It is also assumed that agent $i$ knows a
suitably defined nonlinear map $M_i:\R^n\rightarrow\R^n$ and that all of the $M_i$
share at least one common fixed point.  In general terms, the problem of interest is to develop algorithms, one for each agent,
 which will enable all $m$ agents to iteratively compute a common fixed point of all of the $M_i$.

Motivation for this problem stems, in part, from \cite{lineareqn} which deals with the problem
of devising a distributed algorithm for finding a  solution to the linear equation $Ax = b$,
assuming the equation has at least one solution,
and agent $i$ knows a pair of the matrices
$(A_i^{n_i\times n},b_i^{n_i\times 1})$
where $A = \matt{A_1' &A_2' &\cdots & A_m'}'$ and $b = \matt{b_1'&b_2'&\cdots & b_m'}'$.
Assuming each $A_i$ has linearly independent rows, one local update rule for solving this problem is of the form
$$x_i(t+1) = L_i(z_i(t))$$
where $L_i:\R^n\rightarrow\R^n$ is the affine linear map $x\longmapsto x-A_i'(A_iA_i')^{-1}(A_ix - b_i)$,
$$z_i(t) = \frac{1}{m_i(t)}\sum_{j\in\mathcal{N}_i(t)} x_j(t),$$
and $m_i(t)$ is the number of labels in $\mathcal{N}_i(t)$ (\cite{ACC16.1}).
The map $L_i$ is an example of a
`paracontraction'  with respect to the two norm on $\R^{n}$. More generally,
a continuous nonlinear map $M:\R^n\rightarrow\R^n$  is a {\em paracontraction}
with respect to a given norm $\|\cdot\|$ on $\R^n$,
if $\|M(x)-y\|< \|x-y\|$ for all $x\in\R^n$ satisfying $x \neq M(x)$ and all $y \in \R^n$ satisfying $y = M(y)$
(\cite{paracontract}).
One obvious consequence of this definition is that
$\|M(x)-y\| \leq \|x-y\|$ for all $x\in\R^n$ and all $y \in \R^n$ satisfying $y = M(y)$.
Note that $y=L_i(y)$ if and only if $A_i y = b_i$
and for any such $y$,
$L_i(x) - y = P_i(x-y)$ where $P_i$ is the orthogonal projection matrix $P_i =I-A_i'(A_iA_i')^{-1}A_i$.
 Since the induced 2-norm of $P_i$ is $1$,  $\|P_i(x-y)\|_2\leq \|x-y\|_2,\;\forall x,y\in\R^n$ so
$\|L_i(x)-y\|_2\leq \|x-y\|_2,\;\forall x,y\in\R^n$.
Moreover for any $y$ satisfying $L_i(y)=y$, the inequality $x\neq L_i(x)$ is equivalent
to $x-y\notin\ker A_i$ and $\ker A_i = {\rm image }\;P_i$ so
 $x-y\notin{\rm image }\;P_i$ whenever $x\neq L_i(x)$ and $y\in{\rm image }\;P_i$.
  But for such $x$ and $y$, $\|P_i(x-y)\|_2<\|x-y\|_2$
 so $\|L(x)-y\|_2 <\|x-y\|_2$.  Clearly $L_i$ is a paracontraction  as claimed.

There are many other examples of paracontractions discussed in the literature. Some can be found in \cite{paracontract}
and \cite{Byrne2007}. 
Here are several others.
\begin{enumerate}
    \item The orthogonal projector $x\mapsto \arg\min_{y \in \mathcal{C}} \|x - y\|_2$ associated with a nonempty
        closed convex set $\mathcal{C}$.
This been used for a number of applications including the constrained consensus problem in \cite{Nedic2010}. 
        The fixed points of this map are vectors in $\mathcal{C}$. (\cite{paracontract})
    \item The gradient map $x\longmapsto x - \alpha \nabla f(x)$ where $f :  \R^n \longrightarrow
     \R$ is convex and differentiable, $\nabla f$ is Lipschitz
     continuous with parameter $\lambda>0$, and $\alpha$ is a constant satisfying $0
     < \alpha < \frac{2}{\lambda}$. 
        The fixed points of this map are vectors in $\R^n$ which minimize $f$.
    \item The proximal map associated with a closed proper convex function $f : \R^n \to (-\infty,\infty]$.
        The fixed points of this map are vectors in $\R^n$ which minimize $f$. See \cite{Eckstein1992} as well as \cite{Parikh2014}.
\end{enumerate}
Paracontractions are also discussed in~\cite{Xiao2006} and \cite{Wu2007}.
 What is especially important about paracontractions, whether they are linear or not, is the following well-known theorem published in \cite{paracontract}.
\begin{thm}
Let $M_1,M_2,\ldots,M_m$, be a finite set of $m$ paracontractions with respect to any given norm on $\R^n$.
Suppose that all of the paracontractions share at least one common fixed point.
Let $\sigma(t),\ t \in \{1, 2, \ldots\}$ be an infinite sequence of integers from the set $\{1, 2, \ldots, m\}$
with the property that each integer in $\{1, 2, \ldots, m\}$ occurs in the sequence infinitely often.
Then the state $x(t)$ of the iteration
$$x(t+1) = M_{\sigma(t)}(x(t)),\;\;\;\;t\in\{1,2,\ldots\}$$
converges to a common fixed point of the $m$ paracontractions.\label{eisner}\end{thm}
In the sequel we will use this result to establish the convergence of a family of distributed paracontracting iterations.

  \section{The Problem}

The specific problem to which this paper is addressed is this.
Let  $M_1,M_2,\ldots, M_m$ be a  set of $m$ paracontractions
   with respect to the standard $p$-norm $\|\cdot\|$ on $\R^n$ where $p$ is a constant satisfying $1<p<\infty$.
    Suppose that all of the paracontractions   share at least one common
    fixed point. Find conditions on the time-varying neighbor graph $\mathbb{N}(t)$ so that the states of all $m$ iterations
\eq{x_i(t+1) = M_i\left (\frac{1}{m_i(t)}\sum_{j\in\mathcal{N}_i(t)} x_j(t)\right),\;\;i\in\mathbf{m},\;t\geq 0\label{iter}}
converge to a common fixed point of the $M_i$ where $\mathbf{m} \dfb \{1, \ldots, m\}$ and $\mathcal{N}_i(t)$ is the set of labels of
those agents which are neighbors of agent  $i$ at time $t$. The main result of this paper is as follows.

\begin{thm} If each of the neighbor graphs in the sequence $\mathbb{N}(1),\mathbb{N}(2),\ldots $ is strongly connected
and the paracontractions $M_1,M_2,\ldots, M_m$ share at least one common fixed point, then the states $x_i(t)$
 of the $m$ iterations defined by  \rep{iter}, all converge  to a common fixed point of the $M_i$
 as $t\rightarrow\infty$.
\label{main}\end{thm}
The remainder of this paper is devoted to a proof of this theorem.

\section{Analysis}

 To proceed, let us note that the family of $m$ iterations given by \rep{iter} can be written
 as a single iteration of the form
 \eq{x(t+1) = M((F(t)\otimes I)x(t)),\;\;\;\;t\geq 0\label{big}}
where  for any set of  vectors $\{x_i\in\R^n,\;i\in\mathbf{m}\}$, $x\in\R^{mn}$ is the stacked vector
\eq{x= \matt{x_1 \\ x_2 \\ \vdots \\ x_m}\label{stack}}
$M:\R^{mn}\rightarrow\R^{mn}$ is the map
$$x \longmapsto  \matt{M_1(x_1) \\ M_2(x_2) \\ \vdots \\ M_m(x_m)},$$
$F(t)$ is the  $m\times m$ flocking matrix\footnote{By the {\em flocking matrix} of a
  neighbor graph $\mathbb{N}$
  is meant that stochastic matrix $F=D^{-1}A'$ where $A$ is the adjacency matrix of $\mathbb{N}$,  $D$
  is a diagonal matrix whose $i$th diagonal entry is the in-degree of
   vertex $i$ in $\mathbb{N}$, and prime denotes transpose.} determined by $\mathbb{N}(t)$, $I$ is the $n\times n$ identity matrix, and
 $F(t)\otimes I$ is the Kronecker product of $F(t)$ with $I$.

It will be convenient to introduce the  ``average'' vectors
\eq{z_i(t) = \frac{1}{m_i(t)}\sum_{j\in\mathcal{N}_i(t)} x_j(t),\;\;\;\;i\in\mathbf m,\;\;t\geq 0\label{average}}
in which case  the stacked vector
$$z(t) = \matt{z_1'(t) &z_2'(t)&\cdots &z_m'(t)}'$$
satisfies
\eq{z(t) = (F(t)\otimes I)x(t),\;\;\;\;t\geq 0\label{big22}}
and consequently
\eq{z(t+1) = (F(t+1)\otimes I)M(z(t)),\;\;\;\;t\geq 0\label{big2}}
because of \rep{big}.
It is clear that convergence of all of the $x_i$ to a single point in $\R^n$ implies
 convergence of all of the
$z_i$ to the same point. On the other hand, if all of the $z_i$ converge to a single point which is,
 in addition, a common fixed point of the $M_i,\;i\in\mathbf{m}$, then because
the $M_i$ are continuous and $x_i(t+1) = M_i(z_i(t)),\;t\geq 0$, all of the $x_i$ must  converge
 to the same fixed point. In other words,
 convergence of all of the $z_i$ to a common fixed point of the
 $M_i,\;i\in\mathbf{m}$, is equivalent to convergence of all of the $x_i$ to the same fixed point.
Thus to prove Theorem~\ref{main} it is enough to show that
  if all of the  $\mathbb{N}(t)$ are strongly connected,  the $z_i(t)$ all converge to a
common fixed point $y^*$ of the $M_i,\;i\in\mathbf{m}$.

It is obvious from \rep{big2} that for any positive integer $q$,
\eq{z(q) = ((F(q)\otimes I)M\circ\cdots \circ (F(1)\otimes I)M)(z(0))\label{zform}}
Prompted by this we will study the properties of  maps from $\R^{mn}$ to $\R^{mn}$ which are of the form
$x\longmapsto ((S(q)\otimes I)M\circ\cdots \circ (S(1)\otimes I)M)(x)$ where $q$ is a positive integer,
and $S(t)$, $t\in\mathbf{q} \dfb \{1,2,\ldots,q\}$
is a family of $q$ stochastic matrices $S(t) = \matt{s_{ij}(t)}_{m\times m}$.
We will show that
  under suitable conditions, such maps  are paracontractions with respect to the {\em mixed vector norm}
 $\|\cdot\|_{p,\infty}$ on $\R^{mn}$ where $p$ is a value satisfying $1 < p < \infty$.
and for stacked vectors $x$ of the form shown in \rep{stack},  $$ \|x\|_{p,\infty} = \max_{i\in\mathbf{m}} \|x_i\|_p $$
Here $\|\cdot\|_p$ is the standard $p$ norm on $\R^n$.
The main technical result of this paper is as follows.

\begin{thm}
    Let $M_i,\;i\in\mathbf{m}$ be a set of $m>1$ paracontractions with respect to the standard $p$ norm
 $\|\cdot\|_p$ on $\R^n$ where $p$ is a constant satisfying $1<p<\infty$. Let $S(1),S(2),\ldots, S(q)$ be a set of $q\geq 1$
$m\times m$ stochastic matrices.  If the $M_i,\;i\in\mathbf{m}$ have a common fixed point and
 the matrix product $S(q)S(q-1)\cdots S(1)$ is positive, then the composed map
$\R^{mn}\rightarrow \R^{mn}$, $x\longmapsto ((S(q)\otimes I)M\circ\cdots \circ (S(1)\otimes I)M)(x)$
\begin{enumerate}
\item \label{jel1} is a paracontraction
 with respect to the mixed vector norm $\|\cdot\|_{p,\infty}$.
\item \label{jel2} has as its set of fixed points  all stacked vectors of the form $\matt{y' & y'&\cdots & y'}'$
 where $y$ is a common fixed point of the $M_i,\;i \in \mathbf{m}$.
 \end{enumerate}
\label{tech}\end{thm}
This theorem will be proved later in this section.
In order to prove Theorem~\ref{main}, we will need the following lemma.
\begin{lem}\label{S-qne}
Let $S = \matt{s_{ij}}_{m \times m}$ be a stochastic matrix.
Then
\eq{\| (S \otimes I)x - \bar y \|_{p,\infty} \le \| x - \bar y \|_{p,\infty}}
for any $x \in \R^{mn}$ and $\bar y \in \R^{mn}$ of the form $\bar y = \matt{y' & \cdots & y'}'$.
\end{lem}
\noindent{\bf Proof of Lemma~\ref{S-qne}.}
For each $i \in \mathbf{m}$,
$$
\begin{array}{rl}
\left\Vert \sum_{j \in \mathbf{m}} s_{ij} x_j(t) - y \right\Vert_p
&\le \sum_{j \in \mathbf{m}} s_{ij} \| x_j - y \|_p \\
&\le \left( \sum_{j \in \mathbf{m}} s_{ij} \right)\max_{j \in \mathbf m} \| x_j - y \|_p
\end{array}
$$
by the triangle inequality and that fact that $0 \le s_{ij} \le 1$ for each $i \in \mathbf{m}$ and $j \in \mathbf{m}$.
But since $S$ is stochastic, $\sum_{j \in \mathbf{m}} s_{ij} = 1$ as well.
Thus for each $i \in \mathbf{m}$,
$$
\left\Vert \sum_{j \in \mathbf{m}} s_{ij} x_j(t) - y \right\Vert_p
\le \max_{j \in \mathbf m} \| x_j - y \|_p.
$$
Since this holds for each $i \in \mathbf{m}$,
$$
\max_{i \in \mathbf{m}} \left\Vert \sum_{j \in \mathbf{m}} s_{ij} x_j(t) - y \right\Vert_p
\le \max_{j \in \mathbf m} \| x_j - y \|_p.
$$
Therefore,
$$\|(S \otimes I)x - \bar y \|_{p,\infty} \le \| x - \bar y \|_{p,\infty}.$$
$\qed$

\noindent{\bf Proof of Theorem~\ref{main}.}
All neighbor graphs in the sequence $\mathbb{N}(1), \mathbb{N}(2), \ldots $
 have self arcs at all vertices because each agent is assumed to be a neighbor of itself. It is known that the composition of
$n-1$ such graphs must be complete because each of the graphs in the sequence is, by assumption, strongly connected
\{c.f., Proposition 4 of \cite{reachingp1}\}.
This means that the product of any $q \dfb n-1$ flocking matrices $F(t)$  must be  a positive matrix.
Thus for each $i \ge 1$, the matrix $F(iq) \cdots F(1+(i-1)q)$ is positive.
From \rep{big2}, it follows that
\begin{equation}\begin{array}{rl}
z(iq)
&= ((F((iq) \otimes I) M \circ \cdots \\
& \quad \cdots \circ (F(1 + (i - 1)q) \circ I)M)(z((i - 1)q))
\end{array}\label{subsample-iter}\end{equation}
for each $i \ge 1$.
It follows from  Assertion~\ref{jel1} of Theorem~\ref{tech}
that the maps  $x\longmapsto ((F(iq)\otimes I)M\circ\cdots \circ (F(1 +(i-1)q)\otimes I)M)(x),\ i\geq 1$  are all
paracontractions with respect to the mixed vector norm $\|\cdot \|_{p,\infty}$. 
Moreover there are only finitely many such maps
because there are only a finite number of $n\times n$ flocking matrices.
Furthermore it is clear from Assertion~\ref{jel2} of Theorem~\ref{tech},
that any fixed point $\bar{y}$ common to  these maps  is of the  form $\bar{y} = \matt{y' & y'&\cdots & y'}'$
where $y$ is a common fixed point of the $M_i,\;i \in \mathbf{m}$.
It is clear from Theorem~\ref{eisner} and \rep{subsample-iter}
that $z(iq),\ i \ge 0$ must converge to such a fixed point $\bar{y}$.

By Lemma \ref{S-qne}, $\| z(t+1) - \bar y \|_{p,\infty}
= \| (S(t+1) \otimes I) M (z(t)) - \bar y \|_{p,\infty}
\le \| M (z(t)) - \bar y \|_{p,\infty}$
for any $t \ge 0$.
But since each $M_i,\ i \in \mathbf{m}$ is paracontracting and $y$ is a common fixed point,
$\| M(z(t)) - \bar y \|_{p,\infty}
= \max_{i \in \mathbf m} \| M_i(z_i(t)) - y \|_p
\le \max_{i \in \mathbf m} \| z_i(t) - y \|_p
= \| z(t) - \bar y \|_{p,\infty}
$
Thus,
\eq{\| z(t+1) - \bar y \|_{p,\infty} \le \|z(t) - \bar y\|_{p,\infty}, \quad t \ge 0.\label{z-fejer-monotone}}
From this and the fact that $z(iq),\ i \ge 0$ converges to $\bar y$,
it is also true that $z(t),\ t \ge 0$ must also converge to $\bar{y}$.
Consequently, each $z_i(t)$ must converge to the same common fixed point $y$.
It follows that each $x_i(t)$ must converge to $y$ as well.
$\qed$

In the sequel we develop the technical results needed to prove Theorem~\ref{tech}.
In doing this we will make use of the matrix
$\Phi(t,\tau) = \matt{\phi_{ij}(t,\tau)}_{m \times m}$ which we define as
$\Phi(t,\tau) = S(t)S(t-1)\cdots S(\tau+1)$ for $0 \leq \tau < t\leq q$
and 
$\Phi(t,t) = I$ for $0 \le t \le q$.
Note that
 $S(t)\Phi(t-1,\tau) = \Phi(t,\tau) = \Phi(t,\tau+1)S(\tau+1),\;0\leq \tau<t\leq q$.
For each $i\in\mathbf{m}$, let $v_i(0)\in\R^n$ be an arbitrary but fixed vector, and define
\eq{v_i(t+1) = \sum_{j\in\mathbf{m}}s_{ij}(t+1)M_j(v_j(t)),\quad 0 \le t < q.\label{v-defn}}
 We shall need the following lemmas.

\begin{lem}\label{lem-ineq}
Let $y^*$ be a common fixed point of the $M_i,\;i\in\mathbf{m}$.
For each $i\in\mathbf{m}$
\eq{\|v_i(t)-y^*\| \leq \sum_{j\in\mathbf{m}}\phi_{ij}(t,\tau) \|v_j(\tau) - y^*\|, \label{br1}}
for $0\leq \tau\leq t\leq q$.
\end{lem}
\noindent{\bf Proof of Lemma~\ref{lem-ineq}.}
Fix $0 \le \tau \le q$.
If $t = \tau$, then \rep{br1} holds for each $i\in\mathbf{m}$ since $\phi_{ij}(t,\tau) = 1$ whenever $i=j$ and $\phi_{ij}(t,\tau) = 0$ whenever $i \neq j$.

Suppose $t > \tau$ and \rep{br1} holds for some $t = \mu$ satisfying $\tau \le \mu < q$,
\eq{\|v_i(\mu)-y^*\| \leq \sum_{j\in\mathbf{m}}\phi_{ij}(\mu,\tau) \|v_j(\tau) - y^*\|,\quad i \in \mathbf{m}. \label{lem1-1}}
From \rep{v-defn} and the triangle inequality it follows that
$\|v_i(\mu+1)-y^*\|\leq \sum_{j\in\mathbf{m}}s_{ij}(\mu+1)\|M_j(v_j(\mu))-y^*\|$.
But the $M_i$ are paracontractions, so
\eq{\|v_i(\mu+1)-y^*\| \leq \sum_{j\in\mathbf{m}}s_{ij}(\mu+1)\|v_j(\mu)-y^*\|,\;i \in \mathbf{m}. \label{lem1-2}}
From \rep{lem1-1} and \rep{lem1-2}, it follows that
$$\begin{array}{rl}
    \|v_i(\mu+1)-y^*\| \!\!
    & \leq \! \sum_{k\in\mathbf{m}} s_{ik}(\mu+1) \! \sum_{j\in\mathbf{m}} \! \phi_{kj}(\mu,\tau)\|v_j(\tau)-y^*\| \\
    & = \! \sum_{j\in\mathbf{m}} \sum_{k\in\mathbf{m}} \! s_{ik}(\mu+1) \phi_{kj}(\mu,\tau) \|v_j(\tau)-y^*\|
\end{array} $$
for each $i \in \mathbf m$.
But $\phi_{ij}(\mu+1,\tau) = \sum_{k\in\mathbf{m}} s_{ik}(\mu+1) \phi_{kj}(\mu, \tau)$ by the definition of $\Phi$,
so
$$\|v_i(\mu+1)-y^*\|\leq\sum_{j\in\mathbf{m}}\phi_{ij}(\mu+1,\tau) \|v_j(\tau)-y^*\|, \ i \in \mathbf m$$
which shows that \rep{br1} holds for $t=\mu+1$.
By induction, \rep{br1} holds for any $t$ satisfying $\tau < t \le q$.
Since $\tau$ was initially fixed, \rep{br1} holds for any $0 \le \tau \le t \le q$.
$\qed$

\begin{lem}\label{lem-strict-ineq}
Let $y^*$ be a common fixed point of the $M_i,\; i\in\mathbf{m}$.
Then for each $i\in\mathbf{m}$,
\eq{\|v_i(q)-y^*\|\leq \sum_{j\in\mathbf{m}}\phi_{ij}(q,0)\|v_j(0)-y^*\|\label{gum1}}
and the following statements are true.
\begin{enumerate} 
\item\label{gum2}
If there is a $t$ satisfying $0 \le t < q$ and a $j\in\mathbf{m}$
for which $\phi_{ij}(q,t)>0$ and $M_j(v_j(t))\neq v_j(t)$, then
\eq{\|v_i(q)-y^*\|< \sum_{p\in\mathbf{m}}\phi_{ip}(q,0)\|v_p(0)-y^*\|.\label{pup}}
\item\label{gum3}
If for every $t$ satisfying $0 \le t < q$ and $j\in\mathbf{m}$
it is true that $M_j(v_j(t)) = v_j(t)$
whenever $\phi_{ij}(q,t)>0$,
then
\eq{v_i(q)= \sum_{p\in\mathbf{m}}\phi_{ip}(q,0)v_p(0).\label{boop}}
\end{enumerate}
\end{lem}

\noindent{\bf Proof of Lemma~\ref{lem-strict-ineq}.}
Fix $i\in\mathbf{m}$.
Observe that by setting $t=q$ and $\tau = 0$ in \rep{br1}, one obtains \rep{gum1}.
To prove Assertion~\ref{gum2}, fix $t$ to satisfy $0 \le t < q$ and $j\in\mathbf{m}$ and suppose that $\phi_{ij}(q,0)>0$ and $M_j(v_j(t)) \neq v_j(t)$.
The latter implies that
\eq{\|M_j(v_j(t))-y^*\|<\|v_j(t)-y^*\|\label{loom}}
since $M_j$ is a paracontraction.
From \rep{v-defn} and the triangle inequality
\eq{\|v_p(t+1)-y^*\|\! \leq \! \sum_{k\in\mathbf{m}} \! s_{pk}(t+1)\|M_k(v_k(t))-y^*\|,\;p\in \mathbf{m}\label{ppp}}
By \rep{br1}
$$\|v_i(q)-y^*\| \leq \sum_{p\in\mathbf{m}}\phi_{ip}(q,t+1) \|v_p(t+1) - y^*\|.$$
This and \rep{ppp} imply that
$$\|v_i(q)-y^*\| \leq \sum_{k\in\mathbf{m}}\sum_{p\in\mathbf{m}} \phi_{ip}(q,t+1)  s_{pk}(t+1)\|M_k(v_k(t))-y^*\|.$$
But $\phi_{ik}(q,t) = \sum_{p \in \mathbf m} \phi_{ip}(q,t+1) s_{pk}(t+1)$, so
\eq{\|v_i(q)-y^*\| \leq  \sum_{k\in\mathbf{m}}\phi_{ik}(q,t)\|M_k(v_k(t))-y^*\|.\label{oop}}

Note that \rep{oop} can be written as
\eq{\begin{array}{rl}
    \|v_i(q)-y^*\| \leq
    &\phi_{ij}(q,t)\|M_j(v_j(t))-y^*\| \\
    & + \sum_{\substack{k \in \mathbf{m} \\ k \neq j}} \phi_{ik}(q,t)\|M_k(v_k(t))-y^*\|.
\end{array}\label{oop2}}
By hypothesis, $\phi_{ij}(q,t)>0$.
Moreover $\|M_k(v_k(t))-y^*\|\leq\|v_k(t)-y^*\|,\;k\in\mathbf{m},$
because each $M_k$ is a paracontraction.
From this, \rep{loom} and \rep{oop2} it follows that
\eq{\|v_i(q)-y^*\| < \sum_{k\in\mathbf{m}}\phi_{ik}(q,t)\|v_k(t)-y^*\|.\label{popp}}
By \rep{br1}
$$\|v_k(t)-y^*\| \leq \sum_{p\in\mathbf{m}}\phi_{kp}(t,0)\|v_p(0) - y^*\|,\;k\in\mathbf{m}.$$
From this and \rep{popp} it follows that
$$\|v_i(q)-y^*\| < \sum_{p\in\mathbf{m}} \sum_{k\in\mathbf{m}}\phi_{ik} (q,t)\phi_{kp}(t,0)\|v_p(0) - y^*\|.$$
But $\phi_{ip}(q,0) = \sum_{k\in\mathbf{q}}\phi_{ik}(q,t)\phi_{kp}(t,0)$  so
$$\|v_i(q)-y^*\| < \sum_{p\in\mathbf{m}} \phi_{ip} (q,0)\|v_p(0) - y^*\|.$$
Therefore, Assertion~\ref{gum2} is true.

We turn to the proof of Assertion~\ref{gum3}.
We claim that for all $t$ satisfying $0 \le t < q$ and all $j\in\mathbf{m}$,
\eq{\phi_{ij}(q,t)M_j(v_j(t)) = \phi_{ij}(q,t)v_j(t).\label{viv}}
This is obviously true if $\phi_{ij}(q,t) = 0$.
But \rep{viv} also holds if $\phi_{ij}(q,t)\neq 0$, because of the hypothesis of Assertion~\ref{gum3}.

Since $\phi_{ip}(t,t) = 1$ whenever $i=p$ and $\phi_{ip}(t,t) = 0$ whenever $i \neq p$,
\eq{v_i(q) = \sum_{p\in\mathbf{m}}\phi_{ip}(q,t)v_p(t),\;\;\;i\in\mathbf{m}.\label{viv3y}}
holds for $t=q$.
Suppose \rep{viv3y} is true for some $t$ satisfying $0 < t \le q$.

From \rep{v-defn}
$$v_p(t) = \sum_{k\in\mathbf{m}}s_{pk}(t)M_k(v_k(t-1)),\;\;\;p\in\mathbf{m}.$$
From this and \rep{viv3y} there follows
$$v_i(q) = \sum_{k\in\mathbf{m}}\sum_{p\in\mathbf{m}}\phi_{ip}(q,t)s_{pk}(t)M_k(v_k(t-1)),\;\;\;i\in\mathbf{m}.$$
But $\phi_{ik}(q,t-1) = \sum_{p\in\mathbf{m}}\phi_{ip}(q,t)s_{pk}(t)$ so
$$v_i(q) = \sum_{k\in\mathbf{m}}\phi_{ik}(q,t-1)M_k(v_k(t-1)),\;\;\;i\in\mathbf{m}.$$
This and \rep{viv} imply that
$$v_i(q) = \sum_{k\in\mathbf{m}}\phi_{ik}(q,t-1)v_k(t-1),\;\;\;i\in\mathbf{m}.$$
Hence by induction, \rep{viv3y} holds for all $t$ satisfying $0 \le t \le q$.
Setting $t=0$ yields \rep{boop}
so Assertion~\ref{gum3} of the lemma is true. $\qed $

In the sequel, $\mathcal{F}(Q)$ denotes the set of fixed points of the map
$Q:\R^n\rightarrow\R^n$; i.e., $\mathcal{F}(Q) = \{x:Q(x)=x\}$.
Additionally, $\mathcal C \subset \mathbb R^{mn}$ denotes the \textit{consensus set},
$\mathcal{C} = \{ \matt{x_1' & x_2' & \cdots & x_m'}' : x_i = x_j,\  i,j \in \mathbf m \}$

\begin{lem} If the matrix product $S(q)S(q-1)\cdots S(1)$ has a strongly connected graph, then
$$\mathcal{F}((S(q)\otimes I)M \circ \cdots \circ (S(1)\otimes I)M) = \mathcal{F}(M) \cap \mathcal{C}$$
where $\circ$ denotes composition.
\label{class}\end{lem}

\noindent{\bf Proof of Lemma~\ref{class}.}
Let $x\in \mathcal{F}(M)\cap\mathcal{C}$.
Therefore $x\in\mathcal{C}$ and all of the subvectors  $x_i$ of $x =\matt{x_1'& x_2'&\cdots &x_m'}'$ must be equal.
This in turn implies that $(S(t)\otimes I)x = x,\;t\in\mathbf{q}$ since each $S(t)$ is a stochastic matrix.
Since $x\in\mathcal{F}(M)$, $M(x) = x$.
Thus $(S(t)\otimes I)M(x) = x,\;t\in\mathbf{q}$ so $((S(q)\otimes I)M\circ \cdots \circ (S(1)\otimes I)M)(x) = x$.
Hence $x\in \mathcal{F}((S(q)\otimes I)M\circ\cdots \circ (S(1)\otimes I)M )$
and thus $\mathcal{F}(M)\cap\mathcal{C}\subset \mathcal{F}((S(q)\otimes I)M\circ\cdots \circ (S(1)\otimes I)M)$.

For the reverse inclusion, let $x\in \mathcal{F}((S(q)\otimes I)M\circ\cdots \circ (S(1)\otimes I)M )$.
Set $v(0) = x$ and let $v(t)=\matt{v_1'(t) &v_2'(t) &\cdots & v_m'(t)}',\;0 \le t \le q$,
where $v_i(0) = x_i,\;i\in\mathbf{m}$ and for $t \in \mathbf{q}$, each $v_i(t)$ is as defined in \rep{v-defn}.
Then $v(q) = v(0) = x$.
Let $y^*$ be a common fixed point of the $M_i,\;i\in \mathbf{m}$.
In view of \rep{gum1},
$$\|v_i(q)-y^*\|\leq \sum_{j\in\mathbf{m}}\phi_{ij}(q,0)\|v_j(0)-y^*\|,\;i\in\mathbf{m}.$$
Thus
$w\leq \Phi(q,0)w$  where $\|v_i(q) -y^*\|$ is the $i$th component of the $n$-vector $w$ and $\leq$ here means component-wise.
Since $\Phi(q,0) = S(q)S(q-1)\cdots S(1)$ has a strongly connected graph, $\Phi(q,0)$ is irreducible.
It follows that  $w = \Phi(q,0)w$ \{c.f., page 530 of \cite{Horn2013}\}.
By the Perron-Frobenius Theorem, all components of $w$ must be the same so all $\|v_i(q)-y^*\|, \; i\in\mathbf{m}$ must have the same value.

Suppose that for some $t$ satisfying $0 \le t < q$ and $i,j\in\mathbf{m}$, $\phi_{ij}(q,t)> 0$ and $M_j(v_j(t))\neq v_j(t)$.
 By Assertion~\ref{gum2} of Lemma~\ref{lem-strict-ineq},
  $$\|v_i(q)-y^*\| <\sum_{p\in\mathbf{m}}\phi_{ip}(q,0)\|v_p(0)-y^*\|.$$
  Since $v(q) = v(0)$, it follows that $v_p(0) = v_p(q)$ and therefore,
$$\|v_i(q)-y^*\| <\sum_{p\in\mathbf{m}}\phi_{ip}(q,0)\|v_p(q)-y^*\|.$$
 Thus
$$\|v_i(q)-y^*\| <\sum_{p\in\mathbf{m}}\phi_{ip}(q,0)\|v_{a}(q)-y^*\|$$
 where $a\in\mathbf{m}$ is such that $\|v_{a}(q)-y^*\| =\max_{p\in\mathbf{m}}\|{v_p(q)-y^*}\|$.
Since $\sum_{p\in\mathbf{m}}\phi_{ip}(q,0) = 1$, $\|v_i(q)-y^*\| <\|v_{a}(q)-y^*\|$.
This contradicts the fact that all of the
  $\|v_i(q)-y^*\|, \; i\in\mathbf{m}$ have the same value.
Therefore for every $t$ satisfying $0 \le t < q$ and $j\in\mathbf{m}$,
it must be true that $M_j(v_j(t))= v_j(t)$ whenever $\phi_{ij}(q,t)>0$.

By hypothesis, the graph of $\Phi(q,0)$ is strongly connected so
for each $j\in\mathbf{m}$ there must be a $k\in\mathbf{m}$ such that $\phi_{kj}(q,0)> 0$.
This implies that $v_j(0)\in\mathcal{F}(M_j),\;j\in\mathbf{m}$.
Therefore $x\in\mathcal{F}(M)$.

Additionally, the hypothesis of Assertion~\ref{gum3} in Lemma~\ref{lem-strict-ineq} is satisfied.
Therefore
 $$v_i(q) =  \sum_{p\in\mathbf{m}}\phi_{ip}(q,0)v_p(0),\;\;i\in\mathbf{m}.$$
Thus $v(q) = (S(q)\otimes I)\cdots(S(1)\otimes I) x$.
But $v(q) = v(0) = x$, so $x = ((S(q) \cdots S(1)) \otimes I)x$.
Since $S(q)\cdots S(1)$ is strongly connected, the Perron-Frobenius Theorem ensures that all of the subvectors $x_i$ of $x = \matt{x_1'& x_2'&\cdots &x_m'}'$ must be equal and thus $x \in \mathcal{C}$.
Therefore $\mathcal{F}((S(q)\otimes I)\cdots(S(1)\otimes I)) \subset \mathcal{F}(M) \cap \mathcal{C}$.
$\qed $

\begin{lem}
Suppose  $S_{m\times m}$ is a positive stochastic matrix.
Then for any scalar $p$ satisfying $1<p<\infty$,
$S\otimes I$ is a paracontraction  with respect to the mixed vector norm  $\|\cdot \|_{p,\infty}$.
\label{strict2}\end{lem}

\noindent{\bf Proof of Lemma~\ref{strict2}.}
Because $S$ is positive it has a strongly connected graph.
By the Perron-Frobenius Theorem,
the set of fixed points of the map $x
 \longmapsto (S\otimes I)x$ is all vectors of the form
 $\bar{y} =\matt{y'&y'&\cdots & y'}'$ for $y\in\R^n$. Let $x = \matt{x_1'&x_2'&\cdots & x_m'}'$ be
 any vector in $\R^{mn}$ which is not a fixed point of $S \otimes I$.
 Then there must exist integers $i$ and $j$ such that $x_i\neq x_j$.
Suppose first that $x_i$ is a scalar multiple of $x_j$; i.e. $x_i = \lambda x_j$ for some scalar $\lambda$.
Without loss of generality assume $|\lambda|<1$, so $\|x_i\|_p< \|x_j\|_p$.
Clearly $\|x_i\|_p <\|x\|_{p,\infty}$ and for all $d \in \mathbf m$, $\|x_d\|_p \leq\|x\|_{p,\infty}$
Then for each $k \in \mathbf m$,
$$\begin{array}{rl}
    \left \|\sum_{d\in\mathbf{m}}s_{kd}x_d\right \|_p
    &\leq \sum_{d\in\mathbf{m}}\|s_{kd}x_d\|_p \\
    &= \sum_{d\in\mathbf{m}}s_{kd}\|x_d\|_p
    < \sum_{d\in\mathbf{m}}s_{kd}\|x\|_{p,\infty}.
\end{array}$$
This strict inequality holds because $S$ is positive, which ensures that $s_{ki} > 0$.
But $\sum_{d\in\mathbf{m}}s_{kd} = 1$ because $S$ is stochastic so
\eq{\left \|\sum_{d\in\mathbf{m}}s_{kd}x_d\right \|_p < \|x\|_{p,\infty}, \qquad k\in\mathbf{m}.\label{turk}}

Now suppose that $x_i$ is not a  scalar multiple of $x_j$.
Then for each $k \in \mathbf m$, $s_{ki}x_i$ is not a scalar multiple of $s_{kj}x_j$.
By Minkowski's inequality, $\|s_{ki}x_i+s_{kj}x_j\|_p <\|s_{ki}x_i\|_p+\|s_{kj}x_j\|_p$
since $s_{ki}$ and $s_{kj}$ are both positive.
So
\eq{\|s_{ki}x_i+s_{kj}x_j\|_p < s_{ki}\|x_i\|_p+s_{kj}\|x_j\|_p, \quad k \in \mathbf{m}.\label{min}}
By the triangle inequality,
$$\left\|\sum_{d=1} s_{kd}x_d\right \|_p \leq \|s_{ki}x_i+s_{kj}x_j\|_p +\sum_{\substack{d\in\mathbf{m} \\ d \neq i,j}} \|s_{kd}x_d\|_p.$$
Thus using \rep{min},
$$\begin{array}{rl}
    \left \|\sum_{d\in \mathbf{m}} s_{kd}x_d\right\|_p
    &< \sum_{d\in\mathbf{m}}\|s_{kd}x_d\|_p
= \sum_{d\in\mathbf{m}}s_{kd}\|x_d\|_p \\
&\leq \sum_{d\in\mathbf{m}}s_{kd} \|x\|_{p,\infty}=\|x\|_{p,\infty}
\end{array}$$
so \rep{turk} holds for this case as well.
But
$$\|(S\otimes I)x\|_{p,\infty} = \max_{k\in\mathbf{m}}\left \|\sum_{d \in \mathbf{m}} s_{kd}x_d\right\|_p$$
so
\eq{\|(S\otimes I)x\|_{p,\infty}<\|x\|_{p,\infty}.\label{drop}}
Note that for any vector $\bar{y}\in\R^{mn}$ which is a fixed point of $S\otimes I$, $x-\bar{y}\notin \mathcal{F}(S\otimes I)$ because
$x\notin\mathcal{F}(S\otimes I)$.
Since \rep{drop} holds for all vectors which are not fixed points of $S\otimes I$,
it must be true that $\|(S\otimes I)x - \bar{y}\|_{p,\infty}<\|x-\bar{y}\|_{p,\infty}$ so $S\otimes I$ is
a paracontraction as claimed.
$\qed$

\noindent{\bf Proof of Theorem~\ref{tech}:}
First, note that $\mathcal{C} = \{ \matt{y_1' & \cdots & y_m'}' : y_i = y_j,\  i,j \in \mathbf m \}$
and $\mathcal{F}(M) = \{\matt{y_1' & \cdots & y_m'}' :\; y_i \in \mathcal{F}(M_i),\; i \in \mathbf{m} \}$.
From this and Lemma~\ref{class} it follows that
$$\begin{array}{l}\mathcal{F}((S(q) \otimes I)M \circ \cdots \circ (S(1) \otimes I)M) \\
\qquad \qquad = \{\matt{y' & \cdots & y'}' :\; y \in \bigcap_{i=1}^m \mathcal{F}(M_i) \}\end{array}$$
Thus Assertion~\ref{jel2} of the theorem is true.

Pick $\bar{y}\in\mathcal{F}((S(q)\otimes I)M\circ\cdots \circ (S(1)\otimes I)M)$
 and $x\notin\mathcal{F}((S(q)\otimes I)M\circ\cdots \circ (S(1)\otimes I)M)$.
In view of Lemma~\ref{class}, either $x\notin\mathcal{F}(M)$ or $x\notin\mathcal{C}$;
moreover $\bar{y}\in\mathcal{F}(M)$ and $\bar{y}\in\mathcal{C}$.
Thus, $\bar{y}$ must be of the form $\bar{y} = \matt{y' &y'&\cdots &y'}$ for some vector $y\in\R^n$.
In addition, $y$ must be a common fixed point of the $M_i,\;i\in\mathbf{m}$.

Set $v_i(0) = x_i,\;i\in\mathbf{m}$ where $\matt{x_1' &x_2'&\cdots x_m'}' = x$
and let $v_i(t),\ t \in \mathbf{q}$ be as defined by \rep{v-defn}.
To complete the theorem's proof, it is sufficient to show that
if $v(0) \notin \mathcal{F}(M)$ or $v(0) \notin \mathcal{C}$, then
\eq{\|v_i(q) - y\|_p< \max_{j\in\mathbf{m}}\|v_j(0) - y\|_p,\;\;i \in \mathbf{m}.\label{bub}}
This is sufficient because \rep{bub} implies $\max_{j \in \mathbf{m}} \|v_j(q) - y\|_p < \max_{j \in \mathbf{m}} \|v_j(0) - y\|_p$,
and therefore $\|((S(q)\otimes I)M \circ \cdots (S(1) \otimes I)M)(v(0)) - y\|_{p,\infty} < \|v(0) - y\|_{p,\infty}$.

Fix $i \in \mathbf{m}$.
We claim that if there is a $t$ satisfying $0 \le t < q$ and a $j\in\mathbf{m}$ for which
$\phi_{ij}(q,t) > 0$ and $M_j(v_j(t)) \neq v_j(t)$ then \rep{bub} holds.
To justify this claim
note first  that
$$\sum_{j \in \mathbf{m} }\phi_{ij}(q,0) \|v_j(0) - y\|_p
\! \leq \!\! \left (\sum_{j \in \mathbf{m}} \phi_{ij}(q,0) \! \right ) \! \max_{j\in\mathbf{m}} \|v_j(0) - y\|_p.$$
But $\sum_{j \in \mathbf m} \phi_{ij}(q,0) = 1$ so
\eq{\sum_{j \in \mathbf{m}} \phi_{ij}(q,0) \|v_j(0) - y\|_p
\leq \max_{j\in\mathbf{m}} \|v_j(0) - y\|_p.\label{pioo}}
If there is a $t$ satisfying $0 \le t < q$ and a $j\in\mathbf{m}$
for which $\phi_{ij}(q,t) > 0$ and $M_j(v_j(t)) \neq v_j(t)$,
then by Assertion~\ref{gum2} of Lemma~\ref{lem-strict-ineq}
$$\|v_i(q)-y\| < \sum_{j \in \mathbf{m}} \phi_{ij}(q,0) \|v_j(0) - y\|_p.$$
Since this and \rep{pioo} imply \rep{bub}, the claim is true.


To complete the proof there are two cases to consider, the first being  when  $v(0) \notin \mathcal{F}(M)$.
In this case there is some $j \in \mathbf{m}$ such that $v_j(0) \notin \mathcal{F}(M_j)$.
By hypothesis, $\Phi(q,0) = S(q) \cdots S(1)$ is positive and so $\phi_{ij}(q,0) > 0$.
Therefore with this value of $j$ and $t=0$,
$\phi_{ij}(q,t) > 0$ and $M_j(v_j(t)) \neq v_j(t)$.
Hence \rep{bub} holds in this case.

Now consider the case when $v(0) \notin \mathcal{C}$.
Note that $\Phi(q,0)\otimes I$ is a paracontraction by Lemma~\ref{strict2} and the assumption that
$\Phi(q,0) = S(q) \cdots S(1)$ is a positive matrix.
Clearly
\eq{\|(\Phi(q,0) \otimes I) v(0) - \bar {y}\|_{p, \infty}
 < \|v(0) - \bar{y}\|_{p,\infty}.
    \label{S-strict}}
In other words
\eq{\max_{j\in\mathbf{m}} \left \|\sum_{k \in m} \phi_{jk}(q,0) v_k(0) - y \right \|_p
< \max_{j\in\mathbf{m}} \left \|v_j(0) - y \right \|_p. \label{S_i-strict}}

As noted in the above claim,
if there is a $t$ satisfying $0 \le t < q$ and a $j\in\mathbf{m}$
for which $\phi_{ij}(q,t) > 0$ and $M_j(v_j(t)) \neq v_j(t)$
then \rep{bub} holds.
If on the other hand, there is no $t$ satisfying $0 \le t < q$ and $j\in\mathbf{m}$
for which $\phi_{ij}(q,t) > 0$ and $M_j(v_j(t)) \neq v_j(t)$
then  Assertion~\ref{gum3} of Lemma~\ref{lem-strict-ineq} applies, and so
$$v_i(q)= \sum_{p\in\mathbf{m}}\phi_{ip}(q,0)v_p(0).$$
Therefore
$$\|v_i(q)-y\|_p = \left\|\sum_{p\in\mathbf{m}}\phi_{ip}(q,0)v_p(0) - y\right\|_p$$
Additionally,
$$\left\| \sum_{p\in\mathbf{m}}\phi_{ip}(q,0)v_p(0) -y\right\| \le \max_{j\in\mathbf{m}} \left \|\sum_{k \in m} \phi_{jk}(q,0) v_k(0) - y \right \|_p$$
 Finally, from this and \rep{S_i-strict}, it follows that \rep{bub} is true.
$\qed$

\section{Concluding Remarks}
It is more or less obvious that all of the results of this paper extend painlessly to the case when the
averages appearing in \rep{iter} are replaced with arbitrary convex combinations,
so long as there are only finitely many such convex combinations.
The results can also be extended to the case
when the sequence of neighbor graphs $\mathbb{N}(1),\mathbb{N}(2),\ldots$ is repeatedly jointly strongly connected (\cite{lineareqn}).
It also appears likely that these results can be generalized to the case when the $m$ agents act asynchronously (\cite{async.lin}).
It would be interesting to determine necessary conditions on the neighbor graph sequence which ensure convergence.
These issues will be addressed in future work.

\bibliography{nolcos}             

\end{document}